\documentclass[12pt]{article}
\usepackage{latexsym}
\usepackage{amssymb, amsmath, xspace, lscape,  latexsym}
\newcommand{\ds}{\displaystyle}

\renewcommand{\v}{\textsf{v}}
\newcommand{\dzy}{\frac{\textsf{v}'(z)-\textsf{v}'(y)}{z-y}}
\newcommand{\sps}{s\frac{d}{ds}}

\newcommand{\bea}{\begin{eqnarray}}
\newcommand{\eea}{\end{eqnarray}}
\def\beq#1#2\eeq{
        \begin{equation}
        \label{#1}
            #2
        \end{equation}}

\newcommand{\Ga}{\Gamma}

\newcommand{\al}{\alpha}

\newcommand{\bt}{\beta}

\newcommand{\rme}{{\rm e}}
\newcommand{\dth}{{\delta^2H_n}}

\renewcommand{\hat}{\widehat}
\renewcommand{\tilde}{\widetilde}

\def\btheor#1\etheor{
        \begin{theor}
            #1
        \end{theor}
    }

    \def\bsled#1\esled{
        \begin{sled}
            #1
        \end{sled}   }

\def\btheor#1\elemma{
        \begin{lemma}
            #1
        \end{lemma}
    }

    \def\bsled#1\esled{
        \begin{sled}
            #1
        \end{sled}   }

\newtheorem{theorem}{Theorem}

\newtheorem{lemma}{Lemma}

\def\hm#1{#1\nobreak\discretionary{}{\hbox{\m@th$#1$}}{}}
\def\mi#1{\discretionary{\hbox{\m@th$#1$}}{\hbox{\m@th$#1$}}{}}

\vspace{4ex}
\begin{document}
\small{\title{Painl\'eve III and a singular linear statistics in
Hermitean random matrix ensembles I.}
\author{Yang Chen\\
        Department of Mathematics\\
        Imperial College London\\
        180 Queen's Gates\\
        London SW7 2BZ UK\\
        Alexander Its\\
        Department of Mathematical Sciences\\
        Indiana University-Purdue University Indianapolis\\
        402 N. Blackford Street\\
        Indianapolis, IN 46202-3216
        USA}

\date{\textsf{29-01-2009}}}
\maketitle
\centerline{\bf{Abstract}} In this paper, we study a certain linear
statistics of the unitary Laguerre ensembles, motivated in part by an
integrable quantum field theory at finite temperature. It transpires
that this is equivalent to the characterization of a
sequence of polynomials orthogonal with respect to the weight
$$
w(x)=w(x,s):=x^{\al}\rme^{-x}\rme^{-s/x},
\quad 0\leq x<\infty,\;\;\;\al>0,\;\;\;\;s>0,
$$
namely, the determination of the associated Hankel
determinant and recurrence coefficients. Here $w(x,s)$ is the Laguerre
weight $x^{\al}\:\rme^{-x}$ 'perturbed' by a multiplicative
factor $\rme^{-s/x},$ which induces an infinitely strong zero at the origin.

 For polynomials orthogonal on the unit circle,
a particular example where there are explicit formulas, the weight
of which has infinitely strong zeros, was investigated
by Pollazcek and Szeg\"o many years ago. Such weights are said to be
'singular' or irregular due to the violation of the Szeg\"o condition.

 In our problem, the linear statistics is a sum of the reciprocal
of positive random variables $\{x_j:j=1,..,,n\};$ $\sum_{j=1}^{n}1/x_j.$

 We show that the moment generating function, or the Laplace
transform of the probability density function of this linear
statistics is expressed as the ratio of Hankel determinants and as
an integral of the combination of a particular third Painlev\'e function.

\vfill\eject

\noindent

\setcounter{equation}{0}
\section{Introduction.}

It is a well known fact that the joint probability density of the
eigenvalues, $\{x_j:j=1,...,n\}$ of any Hermitean matrix ensemble is \cite{Mehta}
\bea
p(x_1,...,x_n)dx_1...dx_n=\frac{1}{D_n[w]}
\prod_{1\leq j<k\leq n}(x_j-x_k)^2\prod_{l=1}^{n}w(x_l)dx_l,
\eea
where $D_n$ the normalization constant reads
\bea
D_n[w]=\frac{1}{n!}\int_{\mathbb{R}_{+}^n}\prod_{1\leq j<k\leq n}(x_j-x_k)^2
\prod_{l=1}^nw(x_l)dx_l.
\eea
For the sake of concreteness the domain of integration is
$\mathbb{R}_{+}^n:=[0,\infty)^n.$

 Here
$w(\geq 0)$ is a weight supported on $\mathbb{R}_{+}.$ Furthermore we suppose
$w$ has moments of all order, that is,
\bea
\mu_{j}:=\int_{\mathbb{R}_{+}} x^jw(x)dx,\quad j\in\{0,1,...\},
\eea
exist.

 We shall see that $D_n$ will play a fundamental role in this paper.

It is also well known that the normalization constant defined above
has two more alternative representations; the first of which as the
determinant of the Hankel or moment matrix,
\bea
D_n[w]&=&\det\left(\mu_{j+k}\right)_{j,k=0}^{n-1}\nonumber\\
&=&\prod_{j=0}^{n-1}h_j.
\eea
The quantity $h_j$ in the second equality of (1.4) is square of the
weighted $L^2$ norm of monic polynomials of degree $j,$ orthogonal with
respect to $w;$
\bea
\int_{\mathbb{R}_+} P_j(x)P_k(x)w(x)dx=h_j\delta_{j,k}\;.
\eea
The random variable known as the linear statistics is a sum of a function of
the random variables $\{x_j:j=1,...,n\};$
$$\sum_{j=1}^{n}f(x_j),$$
whose probability density function is determined by the standard formula,
\bea
\mathbb{P}_f(Q)=\frac{1}{n!}\int_{\mathbb{R}_{+}^n}
p(x_1,..,x_n)\delta\left(Q-\sum_{j=1}^{n}f(x_j)\right)
dx_1...dx_n.
\eea
The moment generating function denoted as $M_f(s)$
(assuming $f(x)>0$ for all $x\in\mathbb{R}_+$)
is the Laplace transform of $\mathbb{P}_f(Q)$ with respect to $s$
and has the following form:
\bea
M_f(s)&=&\int_{0}^{\infty}\mathbb{P}_f(Q)\rme^{-sQ}dQ\\
&=&\frac{\det\left(\mu_{j+k}(s)\right)_{j,k=0}^{n-1}}
{\det\left(\mu_{j+k}(0)\right)_{j,k=0}^{n-1}}
=\frac{D_n[w(.,s)]}{D_n[w(.,0)]}
=\frac{\prod_{j=0}^{n-1}h_j(s)}{\prod_{j=0}^{n-1}h_j(0)},
\eea
where
\bea
\mu_{j}(s):=\int_{0}^{\infty} x^jw(x)\rme^{-sf(x)}dx,\;\;\;j\in\{0,1,..\}
\eea
and
\bea
h_j(s)=\int_{0}^{\infty} P_j^2(x)w(x)\rme^{-sf(x)}dx,
\eea
are moments and
the square of the $L^2$ norm of the polynomials $P_j$
orthogonal with respect to
$w\exp(-sf),$ respectively. Therefore we see that the moment generating
function is the ratio of the Hankel determinant generated by the
'perturbed' weight $w\exp(-sf)$ to the corresponding quantity generated
by the 'original' weight $w.$

Because the moments depend on $s,$
the coefficients of the polynomials
$P_j(z)$ also depend on $s,$ however, we shall not display this
dependence most of the time.

Our monic polynomials are normalized so that:
\bea
P_n(z,s)=z^n+\mathsf{p}_1(n,s)\;z^{n-1}+...+P_n(0,s),
\eea
with $P_0(z,s):=1$ and $\mathsf{p}_1(0,s):=0.$

 We note here Heine's multiple integral
representation of $P_n(z),$
\bea
P_n(z)=\frac{1}{n!}\int_{\mathbb{R}_{+}^{n}}
\prod_{j=1}^{n}(z-x_j)p(x_1,...,x_n)dx_1...dx_n.
\eea

In section 2, a description is given for a pair of ladder operators for
smooth weights. These will lead to a linear second order
ordinary differential equation satisfied by $P_n(z)$
and two fundamental compatibility conditions valid for
all $z\in\mathbb{C}\cup{\infty}.$

We denote the compatibility conditions as $(S_1)$ and $(S_2).$

 The compatibility conditions are essentially
a consequence of the recurrence relations;
\bea\label{rec}
zP_n(z)=P_{n+1}(z)+\al_n(s)P_n(z)+\bt_n(s)P_{n-1}(z),
\eea
together with the 'initial' conditions: $P_0(z)=1,$ and $\bt_0P_{-1}(z)=0,$
and the Christoffel-Darboux formula (also a consequence of (1.13)).

In the (1.13),
$$\al_n(s)\in\mathbb{R},\quad n=0,1,..$$
 and
$$
\bt_n(s)=\frac{h_{n}}{h_{n-1}}=\frac{D_{n+1}D_{n-1}}{D_n^2}>0,\quad n=1,2,..
$$
are the recurrence coefficients.

 An easy consequence of (1.11) and (1.13) is that
\bea
\mathsf{p}_1(n,s)-\mathsf{p}_1(n+1,s)=\al_n(s).
\eea
Taking a telescopic sum of (1.14) together with $\textsf{p}_1(0,s)=0,$ implies
\bea
\sum_{j=0}^{n-1}\al_j(s)=-\mathsf{p}_1(n,s).
\eea
We refer the readers to \cite{Szego} for basic facts about orthogonal polynomials.

In our problem,
\bea\label{linstat}
f(x):=\frac{1}{x},\quad 0\leq x<\infty,
\eea
and the unperturbed  weight is given by the equation,
$$
w_{0}(x) = x^{\alpha}e^{-x}, \quad \alpha > 0.
$$
The compatibility conditions $(S_1)$ and $(S_2),$ and a combination of
these, $(S_2'),$ produces a pair non-linear difference equations, satisfied
by the auxiliary quantities $a_n$ and $b_n.$ See (2.16) and (2.17).

The recurrence coefficients, $\al_n$ and $\bt_n$ are ultimately
expressed in terms of $a_n$ and $b_n.$ See (2.9) and (2.14).

The linear statistics (\ref{linstat}) leads to the weight
\bea
w(x,s)=w_0(x)\rme^{-s/x}:=x^{\al}\rme^{-x}\:\rme^{-s/x}, \quad \alpha > 0, \quad s \geq 0.\nonumber
\eea
Such weights arise from a certain problem in mathematical physics:
An integrable quantum field theory at finite temperature \cite{Luky}.

 In the theory
of orthogonal polynomials, the effect of infinitely strong zeros
on the Hankel determinants, recurrence coefficients and polynomials
themselves is of considerable interest.

For orthogonal polynomials with weight $w$
supported on $[-1,1],$ the classical
Szeg\"o theory gives a comprehensive account of the
large $n$ behavior of the recurrence
coefficients and the polynomials
(both outside $[-1,1]$ and on $(-1,1)$) if $w$
is absolutely continuous and satisfies the Szeg\"o condition,
$$
\int_{-1}^{1}\frac{|\ln w(x)|}{\sqrt{1-x^2}}dx<\infty.
$$
See [pp. 296--312,\cite{Szego}], \cite{ge} and \cite{Nevai}
regarding Szeg\"o's theory.

However, there is a class of orthogonal polynomials discovered
by Pollaczek and extended
by Szeg\"o which is in some sense irregular. See
[pp. 393--400,\cite{Szego}] and \cite{Szego1}
about this class of orthogonal polynomials. The Pollaczek-Szeg\"o
weight behaves like
$$
{\rm exp}\left(-\frac{c}{\sqrt{1-x^2}}\right),\quad c>0,
$$
near $\pm 1,$ and consequently just violates the Szeg\"o condition.

 We reproduce here some of the results of \cite{Szego1}
to illustrate the irregularity.

 Associated with the weight,
\bea
w(x;a,b):=\frac{\rme^{(2\theta-\pi)\phi(\theta)}}{\cosh[\pi\phi(\theta)]}
\nonumber
\eea
where $x:=\cos\theta,$ $0<\theta<\pi$ and
$$
\phi(\theta):=\frac{a\cos\theta+b}{2\sin\theta},
\;\;\;a,\;b\in\mathbb{R},\;a\geq|b|,
$$
are the normalized Pollaczek polynomials $\{p_n(x;a,b)\},$
$$
\int_{-1}^{1}[p_n(x;a,b)]^2w(x;a,b)dx=1.
$$
If $x\to 1,$ then
$$
w(x;a,b)\simeq\:2\:\rme^{(a+b)(1-\pi/\theta)},\;\;{\rm as}\:\:\theta\to 0,
$$
which shows that the weight vanishes exponentially at $x=1.$
An easy computation demonstrates the same behavior at $x=-1.$

The large $n$ behavior of $p_n(x;a,b)$ are as follows:
\bea
{\rm (a)}\:\:p_n(1;a,b)&\sim& n^{1/4}\:\rme^{2\sqrt{a+b}\:\sqrt{n}},
\nonumber\\
{\rm (b)}\:\:p_n(x;a,b)&\sim&n^K\:[x+(x^2-1)]^{n},\:\:
K=K(x)=\frac{ax+b}{2\sqrt{x^2-1}},\:x\notin[-1,1],\nonumber\\
{\rm (c)}\:\:p_n(\cos\theta;a,b)&=&A(\theta)\:\:
\cos\left[n\theta-\phi(\theta)\:\ln n+B(\theta)\right]+\epsilon_n(\theta),
\:\:\lim_{n\to\infty}\epsilon_n(\theta)=0,\:\:0<\theta<\pi.\nonumber
\eea
where the $n-$ independent functions $A(\theta)(>0),$ $B(\theta)$ are
analytic in $(0,\pi).$

This is to be contrasted with the large $n$ behavior of the
normalized Jacobi polynomials
$\{p_n(x)\}$ associated with the weight
$$
(1-x)^{\al}(1+x)^{\bt},\quad\al>-1,\;\bt>-1,\;\;x\in[-1,1],
$$
\bea
{\rm (a')}\:\:p_n(1)&\sim&n^{\al+1/2},\nonumber\\
{\rm (b')}\:\:p_n(x)&\sim&\left[x+\sqrt{x^2-1}\right]^n,\:\:x\notin[-1,1]\nonumber\\
{\rm (c')}\:\:p_n(\cos\theta)&=&A_1(\theta)\:
\cos\left[n\theta+B_1(\theta)\right]+\epsilon_n(\theta),\nonumber
\eea
where $A_1(\theta)$ and $B_1(\theta)$ are functions of the same kind.

The symbol $\thicksim$ indicates that
the ratio of the given quantities approaches a non-zero limit,
while $\backsimeq$ indicates that the limit is $1.$ This is the convention
adopted by \cite{Szego1} and will not be used later.

 Our paper is a first step
in study of the Pollaczek-Szeg\"o type orthogonal polynomials supported in
infinite intervals.

 With reference to the Heine formula, we see that,
\bea
(-1)^nP_n(0,s)=\frac{D_n[w(.,s,\al+1)]}{D_n[w(.,s,\al)]}.
\eea
For the unperturbed or Laguerre weight, we have the explicit
determination
\bea
(-1)^nP_n(0,0)=\frac{\Ga(n+\al+1)}{\Ga(\al+1)},
\eea
since,
\bea
D_n[w(.,0,\al)]=\frac{G(n+1)G(n+\al+1)}{G(\al+1)},
\eea
where $G(z)$ is the Barnes $G-$function that satisfies the functional relation
$G(z+1)=\Ga(z)G(z).$

In section 3, by taking derivative with respect to $s$ on the orthogonality relations
we obtain a pair of differential-difference equations or the
Toda equations. Combining the Toda equations and the non-linear
difference equations obtained in section 2 produce a particular Painl\'eve III satisfied by
$\al_n(s),$up to linear shift in $n.$

The $\tau-$ function for this $P_{III}$ turns out to be intimately related to
the Hankel determinant
\bea
D_n(s)=\det
\left(\int_{0}^{\infty}x^{j+k}\;x^\al\rme^{-x-s/x}\right)_{j,k=0}^{n-1}.
\eea
We also express the recurrence coefficients $\al_n$ and $\bt_n$ in terms of the
logarithmic derivative of Hankel determinant
$$
H_n:=s\frac{d}{ds}\ln D_n(s),
$$
and obtain a functional equation involving
$H_n,\:H_n'$ and $H_n''.$ The resulting second
order non-linear ordinary differential equation satisfied by $H_n$, is recognized to be
the Jimbo-Miwa-Okamoto $\sigma$ - form of our $P_{III}.$

In section 4, we show, with the aid of the non-linear difference equations derived in
section 2, another functional equation involving $H_n,\:H_{n+1}$ and $H_{n-1}.$ We call
the resulting non-linear second order difference equation satisfied by $H_n$
the discrete $\sigma$ - form of $P_{III}.$

In sections 5, the Riemann-Hilbert approach to orthogonal polynomials and
the isomonodromy deformation theory of Jimbo and Miwa, are used
to re-derive the $P_{III}$ and thereby identify the auxiliary quantities, $a_n$
and $b_n$, introduced in section 2, with the objects of the
Jimbo-Miwa isomonodromy theory of Painlev\'e equations.

In section 6, we show that the Hankel determinant is the isomonodromy $\tau-$function
in the sense of Jimbo and Miwa, and put into context of the general theory
of integrable systems the identities derived in sections
2 and 3.

As we are studying an example of orthogonal polynomials where the
otherwise classical weight,
$$
w_0(x)=x^{\al}{\rm e}^{-x},\;\;x\in\mathbb{R}_+, \quad \alpha > 0
$$
is perturbed by an infinitely strong zero,
$$
w(x,s):=w_0(x)\:{\rm e}^{-s/x},\;\;s\geq 0,
$$
the natural questions of interest are about the large $n$ behavior of the Hankel determinant,
recurrence coefficients and the orthogonal polynomials. Such investigations will therefore
provide valuable insights into the asymptotic of the associated Painlev\'e
transcendant. These results will be published in a forthcoming paper \cite{ChenIts}.

We want to emphasize that in this paper and in its follow-up  we do
not claim the introduction of  new concepts. Our main aim is to
investigate a concrete important example  of the linear statistics
which leads to a  strong zero at $x=0$ using the known techniques.
In addition, taking this example as a  ``case study'' we show how
the  apparatus  which are used by the two communities - the
orthogonal polynomial community and the integrable system community,
match with each other.

\setcounter{equation}{0}
\section{Ladder operators and non-linear difference equations.}

The pair of ladder operators has been known to various authors. See, for example,
\cite{Bau}, \cite{ron}, \cite{ladder1}, \cite{ladder2}, \cite{ladder3}, \cite{ladder4},
\cite{magnus2}, \cite{CIs}, \cite{FIK4}, \cite{FIK5},
\cite{magnus1} and \cite{Sho}. In fact, Magnus in \cite{magnus2} noted that such
operators were known to Laguerre.

 Because the associated fundamental compatibility conditions
and their use in the derivation of the Painl\'eve transcendant
\cite{chen-pain1},\cite{chen-pain2}, \cite{FIK4}, \cite{F+W}, \cite{magnus2} are perhaps less
well known we summarize these findings in (2.1--2.5), $(S_1),$ $(S_2)$ and $(S_2')$
in a form which we find particular easy to use. We note here that $(S_1),\;(S_2)$
and $(S_2')$ were also known to Magnus \cite{magnus2} and $(S_2)$ also appeared
in \cite{I+W}. See also \cite{F+W}.

 For polynomials orthogonal on the unite circle the analogues ladder operators
can be found in \cite{Is+W}. See \cite{F+W1} for the circular case applicable
to bi-orthogonal polynomials. The compatibility condition in the circular
case can be found \cite{BC} where it was used to obtain in explicit form
of the Toeplitz determinant with the pure Fisher-Hartwig symbol
and the discriminant of the associated orthogonal polynomials.

The compatibility conditions can also be adapted to the situation where the
weight has discontinuities. See \cite{chen-pain1} and \cite{BC1}.

The ladder operators are,
\bea
\left(\frac{d}{dz}+B_n(z)\right)P_n(z)&=&\bt_nA_n(z)P_{n-1}(z)\label{21}\\
\left(\frac{d}{dz}-B_n(z)-\v'(z)\right)P_{n-1}(z)&=&-A_{n-1}(z)P_n(z) \label{22}
\eea
where
\bea
A_n(z)&=&\frac{1}{h_n}\int_{0}^{\infty}\dzy P_n^2(y)w(y)dy\\
B_n(z)&=&\frac{1}{h_{n-1}}\int_{0}^{\infty}\dzy P_{n}(y)P_{n-1}(y)w(y)dy\\
\v(z)&:=&-\ln w(z),
\eea
and the associated fundamental compatibility conditions are,
$$
B_{n+1}(z)+B_n(z)=(z-\al_n)A_n(z)-\v'(z)\eqno(S_1)
$$
$$
1+(z-\al_n)(B_{n+1}(z)-B_n(z))=\bt_{n+1}A_{n+1}(z)-\bt_nA_{n-1}(z),\eqno(S_2)
$$
valid for all $z\in \mathbb{C}\mathsf{U}\{\infty\}.$ See
\cite{chen1} for a recent derivation of the compatibility conditions.
To arrive at the equations (2.3) and (2.4) we have assumed that $w(0)=w(\infty)=0.$
This is certainly the situation for our problem since we have assume that $\al>0$
and $s \geq 0$.

 Combining suitably $(S_1)$ and $(S_2)$ gives an expression
involving $\sum_{j=0}^{n-1}A_j(z),$ $B_n(z)$ and $\textsf{v}'(z)$
from which further insight into recurrence coefficients may be gained.

 The equation $(S_2')$ may be thought of as the first integral from $(S_1)$
and $(S_2).$

 Although $(S_2')$ first appeared in \cite{magnus2} in a slightly different form, we
present here a derivation of a version which we find useful in practice.

  Multiplying $(S_2)$ by $A_n(z)$ we see that
the r.h.s. of the resulting equation is a first order difference, while the l.h.s., with
$(z-\al_n)A_n(z)$ replaced by $B_{n+1}(z)+B_n(z)+\v'(z)$ is a first order
difference plus $A_n(z).$ Taking a telescopic sum,
together with the initial conditions $B_0(z)=A_{-1}(z)=0$, produces the Lemma,
\begin{lemma}
$$
B_n^2(z)+\v'(z)B_n(z)+\sum_{j=0}^{n-1}A_j(z)=\bt_nA_n(z)A_{n-1}(z).\eqno(S_2')
$$
\end{lemma}
If $\v'(z)$ is rational then so are $A_n(z)$ and $B_n(z).$ See (2.3) and (2.4). Furthermore,
eliminating $P_{n-1}(z)$ from (2.1) and (2.2) it is easy to show that $y(z):=P_n(z)$
satisfies the second order linear ordinary differential equation,
\bea
y''(z)-\left(\v'(z)+\frac{A_n'(z)}{A_n(z)}\right)y'(z)+
\left(B'_n(z)-B_n(z)\frac{A_n'(z)}{A_n(z)}+\sum_{j=0}^{n-1}A_j(z)\right)
y(z)=0.
\eea
Note that $(S_2')$ has been used to simplify the coefficient of $y(z)$ in (2.6).

The equation (2.6) can also be found in \cite{Sho}, albeit in a different form.

For the problem at hand, $f(x)=1/x,\;\;\; x\geq 0,$ the weight and associated
quantities are,
\bea
w(x,s)&=&x^{\al}\rme^{-x-s/x},\nonumber\\
\textsf{v}(z)&=&z+s/z-\al\ln z,\quad \v'(z)=1-s/z^2-\al/z\nonumber\\
\dzy&=&\frac{1}{z}\left(\frac{\al}{y}+\frac{s}{y^2}\right)
+\frac{s}{z^2y}.\nonumber
\eea
Using these we have the next Lemma.
\begin{lemma}
The "coefficients" $A_n(z)$ and $B_n(z)$ appearing in the ladder operators are
\bea\label{andef}
A_n(z)&=&\frac{1}{z}+\frac{a_n}{z^2},\\
B_n(z)&=&-\frac{n}{z}+\frac{b_n}{z^2},\\
\quad a_n&:=&\frac{s}{h_n} \int_{0}^{\infty}\frac{P_n^2}{y}wdy,\quad\quad a_n(0)=0,\nonumber\\
\quad b_n&:=&\frac{s}{h_{n-1}}\int_{0}^{\infty}\frac{P_nP_{n-1}}{y}wdy,\quad\quad b_n(0)=0.\nonumber
\eea
\end{lemma}
{\em Proof.\;} From the definitions of $A_n(z)$ and $B_n(z)$ and
with the identities,
\bea
1&=&\frac{1}{h_n}\int_{0}^{\infty}\left(\frac{\al}{y}+\frac{s}{y^2}\right)P_n^2wdy\nonumber\\
-n&=&\frac{1}{h_{n-1}}\int_{0}^{\infty}\left(\frac{\al}{y}+\frac{s}{y^2}\right)P_nP_{n-1}wdy,\nonumber
\eea
obtained by integration by parts, we find (2.7) and (2.8). $\quad\quad\quad\Box$
\\
\\
We see that at this stage there are 4 unknowns, $\al_n,\;\bt_n,\;a_n$ and $b_n.$ In what
follows we will show how $(S_1)$ and $(S_2')$ can be applied to obtain amongst
other things the pair non-linear difference equations involving $a_n$ and $b_n$ mentioned
earlier.

 On equating the residues on both sides of $(S_1),$ we find,
\bea
\al_n&=&2n+1+\al+a_n\label{29}\\
b_{n+1}+b_n&=&s-\al_n\:a_n.\label{210}
\eea
Carrying out a similar calculation with $(S_2'),$ gives,
\bea
\bt_n&=&n(n+\al)+b_n+\sum_{j=0}^{n-1}a_j\label{211}\\
\bt_n(a_n+a_{n-1})&=&ns-(2n+\al)b_n\label{212}\\
b_n^2-sb_n&=&\bt_na_na_{n-1}\label{213}.
\eea

 The upshot of these equations is
that $\al_n$ and $\bt_n$ are entirely determined by $a_n$ and $b_n,$
where $\al_n$ is simply $a_n$ plus $2n+1+\al.$

Eliminating $a_{n-1}$ from (2.12) and (2.13) we have the next Lemma.
\begin{lemma}
\bea
\bt_na_n^2=[ns-(2n+\al)b_n]a_n-(b_n^2-sb_n).\label{214}
\eea
\end{lemma}
Therefore (2.14) expresses $\bt_n$ in terms of $a_n$ and $b_n$, and importantly
bypasses the finite sum in (2.11). Eliminating $\bt_n$ from (2.11) and (2.14),
an expression can be found for
$\sum_{j=0}^{n-1}a_j,$ in terms of $a_n$ and $b_n. $

 We state this in the next Lemma
\begin{lemma}
\bea
\sum_{j=0}^{n-1}a_j=-n(n+\al)-b_n+\frac{ns-(2n+\al)b_n}{a_n}\label{215}
-\frac{b_n^2-sb_n}{a_n^2}.
\eea
\end{lemma}
Note that because $a_n(0)=b_n(0)=0,$
(2.9) and (2.11) reduce to $\al_n(0)=2n+1+\al,$ and
$\bt_n(0)=n(n+\al)$, respectively, which we recognize to be the recurrence coefficients
of the monic Laguerre polynomials.

 In summary, with reference to (2.13) and (2.14), we obtain two non-linear
 difference equations, satisfied by $a_n$ and $b_n,$
 \bea
 b_{n+1}+b_n&=&s-(2n+1+\al+a_n)a_n,\\
 (b_n^2-sb_n)(a_n+a_{n-1})&=&[ns-(2n+\al)b_n]a_na_{n-1},
 \eea
 to be iterated in $n$ with the initial conditions,
 \bea
 a_0(s)&=&\sqrt{s}\frac{K_{\al}(2\sqrt{s})}{K_{\al+1}(2\sqrt{s})},\\
 b_0(s)&=&0,
 \eea
 where $K_{\al}(z)$ is the MacDonald functions of the second kind.

 We call (2.16) and (2.17) together with the initial conditions (2.18) and (2.19)
 the MacDonald's hierarchy. See also \cite{ron} for the treatment of a class
of semi-classical weights.

 Solutions for $a_n$ and $b_n$ that are rational functions of $2\sqrt{s},$
  are found for $\al=p+1/2,\;\:p\in\mathbb{Z},$ since
$$
K_{p+1/2}(z)=K_{-p-1/2}(z)=\sqrt{\frac{\pi}{2z}}\rme^{-z}\sum_{k=0}^{p}
\frac{(p+k)!}{k!(p-k)!(2z)^k}.
$$

\setcounter{equation}{0}
\section{Toda evolution and Painl\'eve III.}

 By taking derivatives with respect to $s$ on the orthogonality relation gives rise to
 Toda type equations.

Because,
\bea
h_n&=&\int_{0}^{\infty}P_n^2wdy,\nonumber
\eea
we have,
\bea
\sps h_n&=&-s\int_{0}^{\infty}\frac{P_n^2}{y}wdy,\nonumber
\eea
hence,
\bea
\sps\ln h_n&=&-a_n\label{31}\\
\sps\ln \bt_n&=&a_{n-1}-a_{n}.
\eea
We also have,
\bea
0&=&\frac{d}{ds}\int_{0}^{\infty}P_nP_{n-1}wdy\\
&=&h_{n-1}\frac{d}{ds}\textsf{p}_1(n)-\int_{0}^{\infty}\frac{P_nP_{n-1}}{y}wdy\\
s\frac{d\textsf{p}_1(n)}{ds}&=&b_n\\
s\frac{d\al_n}{ds}&=&s\frac{d a_n}{ds}\nonumber\\
&=&b_{n}-b_{n+1}\nonumber\\
&=&\bt_{n}-\bt_{n+1}+\al_n.\label{36}
\eea
The equation (3.6) follows from (2.9), (1.14) and (2.11)

There is another identity involving $\sum_{j=0}^{n-1}a_j:$
\bea
\sps\sum_{j=0}^{n-1}a_j=-b_n,
\eea
which is an immediate consequence of a telescopic sum of the second equality of (3.6).

 We now show that the Hankel determinant $D_n,$ is up to scaling transformation
the $\tau-$function of the Toda-equations.
Let
$$
\tilde{D}_n(s):=s^{-n(n+\al)}D_n(s).
$$
 We find by summing (3.1) ,
\bea
\sps\ln D_n=-\sum_{j=0}^{n-1}a_j, \label{38}
\eea
since
$$\sum_{j=0}^{n-1}\ln h_j=\ln D_n.$$

Applying $\sps$ to (3.8) and keeping in mind (3.7), (2.11) and (3.8) gives
\bea
\sps\left(\sps\ln D_n\right)&=&b_n\nonumber\\
&=&\bt_n-n(n+\al)-\sum_{j=0}^{n-1}a_j\nonumber\\
&=&\bt_n-n(n+\al)+\sps\ln D_n.\nonumber
\eea
The last equation simplifies to
\bea
 s^2\frac{d^2}{ds^2}\ln D_n(s)=\frac{D_{n+1}D_{n-1}}{D_n^2}-n(n+\al),\nonumber
 \eea
 since
 $$\bt_n=\frac{D_{n+1}D_{n-1}}{D_n^2}.$$
In terms of $\tilde{D}_n(s)$ we have,
 \bea
 \frac{d^2}{ds^2}\ln\tilde{D}_n(s)
 =\frac{\tilde{D}_{n+1}\tilde{D}_{n-1}} {\tilde{D}_n^2}.
\eea
The equation (3.9) is the Toda molecule equation \cite{Sogo} and
shows that $\tilde{D}_n(s)$ is the corresponding $\tau-$ function of
the Toda equations (3.2) and (3.6).

As the Hankel determinant is now identified with the
$\tau-$function (see section 6 for more on this issue), we may expect the emergence of a
Painlev\'e equation.
In fact, $a_n(s)$ satisfies a particular $P_{III}.$ To see this, we first
investigate the evolution of $a_n$ and $b_n$ as functions of $s.$
\begin{lemma}
For a fixed $n,$ the auxiliary quantities $a_n$ and $b_n$ satisfy the following
coupled Riccatti equations:
\bea
s\frac{d a_n}{ds}&=&2b_n+(2n+1+\al+a_n)a_n-s\\
s\frac{d b_n}{ds}&=&\frac{2}{a_n}(b_n^2-sb_n)+(2n+\al+1)b_n-ns\;.
\eea
\end{lemma}
{\em Proof.\:\:} The equation (3.10) follows from applying $\sps$ to (2.9) together with
the first equality of (3.6) and with (2.10) to replace
$b_{n+1}$ by $s-\al_na_n-b_n$.

 A little bit more work is required to prove
(3.11). First apply $\sps$ to (2.11),
\bea
s\frac{d\bt_n}{ds}&=&s\frac{d b_n}{ds}+\sps\sum_{j=0}^{n-1}a_j\nonumber\\
&=&s\frac{d b_n}{ds}-b_n\nonumber\\
&=&\bt_{n}a_{n-1}-\bt_na_n\nonumber\\
&=&\frac{b_n^2-sb_n}{a_n}-\left[ns-(2n+\al)b_n-\frac{b_n^2-sb_n}{a_n}\right],\nonumber
\eea
where the last three equalities follow from (3.7), (3.2),
(2.13) and (2.14).$\quad\quad\quad\Box$

The next theorem identifies $a_n$ as a particular third Painlev\'e function.
\begin{theorem}
For a fixed $n\in\{0,1,2,..\}$ the auxiliary quantity $a_n$ satisfies
\bea
a_n''=\frac{(a_n')^2}{a_n}-\frac{a_n'}{s}+(2n+1+\al)\frac{a_n^2}{s^2}+
\frac{a_n^3}{s^2}+\frac{\al}{s}-\frac{1}{a_n}, \label{p3}
\eea
with the initial conditions
\bea
a_n(0)=0,\quad a_n'(0)=\frac{1}{\al},\quad\al>0.
\eea
If $a_n(s):=-\textsf{q}(s),$ then $\textsf{q}(s)$ is $P_{III'}(-4(2n+1+\al),-4\al,4,-4),$
following the convention of \cite{oka}.
\end{theorem}

\noindent
{\em Proof.\:\:}Eliminate $b_n$ from (3.10) and (3.11) gives (3.12). The initial
conditions follows from a straightforward computation. $\quad\quad\quad\Box$

\noindent
{\bf Remark I.\:} If $n=0,$ then $b_0=0$ and (3.10) is solved by
$$
{\sqrt s}\:\frac{K_{\al}(2\sqrt{s})}{K_{\al+1}(2\sqrt{s})}
$$
which is (2.18). We observe that the above also solve (3.12) for $n=0.$

\noindent
{\bf Remark II.\:} An alternative form, obtained by,
$$
a_n(s)=:\frac{s}{X_n(s)},
$$
reads,
\bea
X_n''=\frac{(X_n')^2}{X_n}-\frac{X_n'}{s}-\frac{\al X_n^2}{s^2}
-\frac{2n+1+\al}{s}+\frac{X_n^3}{s^2}-\frac{1}{X_n},
\eea
which is a $P_{III'}(-4\al,-4(2n+1+\al),4,-4).$ If the derivatives in (3.14)
were neglected, then $X_n$ solves the quartic,
\bea
X^4-\al X^3-(2n+1+\al)sX-s^2=0.
\eea
We may interpret an appropriate solution (3.15) as the geometric mean of the end points
of the support of a single interval equilibrium density. This appears in a potential theoretic
minimization problem, the detail of which is in a forthcoming paper \cite{ChenIts}. We note
that another $P_{III}$ associated with the Toeplitz determinant
$$
\det[I_{j-k+\nu}(\sqrt{t})]_{0\leq j,k\leq n-1},
$$
was found in \cite{F+W1} (see also \cite{F+W2}). Here $I_{r}(z)$ is the modified Bessel's function of the first kind.
We note that this Toeplitz determinant can be thought of as a Toeplitz-analog of the
Hankel determinant we are studying in this paper.
Also, the  above Toeplitz determinant with $\nu=0$ appeared \cite{T+W} in connection with a certain ensemble
of $n\times n$ unitary matrices and Ulam's problem in combinatorics.

In the next theorem, we display two alternative integral representation of $D_n,$
in terms of $a_n$ and $X_n.$
\begin{theorem}
\bea
\ln\frac{D_n(s)}{D_n(0)}
&=&\int_{0}^{s}\left[\frac{t}{2}-\frac{1}{4}\left(\frac{t}{a_n}-\al\right)^2
-a_n\left(n+\frac{\al}{2}\right)-\frac{a_n^2}{4}+\frac{1}{4}\left(1-\frac{ta_n'}{a_n}\right)^2
\right]\frac{dt}{t}\\
&=&\int_{0}^{s}\left[\frac{t}{2}-\frac{1}{4}(X_n-\al)^2-\left(n+\frac{\al}{2}\right)\frac{t}{X_n}
-\frac{t^2}{4X_n^2}+\frac{t^2X_n'^2}{4X_n^2}\right]\frac{dt}{t}.
\eea
\end{theorem}

\noindent
{\em Proof:\:\:} From (3.8) and (2.15) we see that the logarithmic derivative of $D_n(s)$
is expressed in terms of $a_n$ and $b_n.$ If we use (3.10) to
eliminate $b_n$ in favor of $a_n$ and $a_n'$
from the resulting equation, then (3.16) follows after some
simplification. The equation (3.17) follows from
the substitution $a_n(s)=:s/X_n(s).$ Note that $D_n(0)$ is given by (1.19). $\quad\quad\Box$
\\
\\
Put
\bea \label{Hndef}
H_n:=\sps\ln D_n.
\eea
In section 6  we will show that the Hankel Determinant $D_{n}(s)$ can be identified
with the Jimbo-Miwa $\tau$-function corresponding
to the solution $a_{n}(s)$ of the Painlev\'e III$'$  equation (\ref{p3}). Namely, we
will show that
\begin{equation}\label{tau50}
D_{n}(s) = \mbox{const}\:\tau(s)\: {\rm e}^{\frac{s}{2}}s^{\frac{n(n+\alpha)}{2}}.
\end{equation}
In its turn, this relation yields the following formula for the quantity $H_{n}$,
 \begin{equation}\label{sigma0}
 H_{n} = \sigma(s) + \frac{s}{2} + \frac{n(n+\alpha)}{2},
 \end{equation}
where $\sigma(s) \equiv s\frac{d}{ds}\ln \tau$ is the Jimbo-Miwa-Okamoto  $\sigma$-function
corresponding to the equation $P_{III'}.$  Therefore,  $H_n,\;H_n'$ and $H_n''$
should satisfy a functional equation
$$
f(H_n,H_n',H_n'',n,s,)=0,
$$
known as the Jimbo-Miwa-Okamoto $\sigma$ - form of our $P_{III'}.$

 With $H_n$ defined above, it is easy to
see from (3.7), (3.8) and (2.11) that
\bea
b_n&=&sH_n'\\
\bt_n&=&n(n+\al)+sH_n'-H_n.
\eea
In the next theorem we state the non-linear
second order ordinary differential equation satisfied
by $H_n.$
\begin{theorem}
If
\bea
H_n:=\sps\ln D_n(s),
\eea
then
\bea \label{322}
(sH_n'')^2=[n-(2n+\al)H_n']^2-4[n(n+\al)+sH_n'-H_n]H_n'(H_n'-1).
\eea
\end{theorem}

{\rm Proof.\:\:} First we re-write (2.14) and (3.11) as
\bea
\bt_na_n+\frac{b_n^2-sb_n}{a_n}&=&ns-(2n+\al)b_n\\
\frac{2}{a_n}(b_n^2-sb_n)&=&sb_n'-b_n+ns-(2n+\al)b_n
\eea
respectively. Eliminate $a_n$ from (3.25) and (3.26) produces
\bea
\left(sb_n'-b_n\right)^2=[ns-(2n+\al)b_n]^2-4\bt_n(b_n^2-sb_n).
\eea
The equation (3.24) follows by substituting $b_n$ and $\bt_n$ from (3.21) and
(3.22) into (3.27).
$\quad\quad\Box$

Hence the recurrence coefficients $\al_n$ and $\bt_n,$ of the orthogonal
polynomials associated with our weight,
$$
w(x,s)=x^{\al}\rme^{-x-s/x},\quad x\in[0,\infty),\quad \al>0,\quad s>0,
$$
are expressed in terms of $H_n,\;H_n'$ and $H_n'',$ as follows:
\bea
\al_n&=&2n+1+\al+\frac{2s(H_n')^2-sH_n'}{sH_n''+n-(2n+\al)H_n'}\\
\bt_n&=&n(n+\al)+sH_n'-H_n,
\eea
and that $H_n$ itself satisfies (3.24).

\vskip .2in
\noindent
{\bf Remark III.\:} With the identification (\ref{sigma0}) of the function $H_{n}$ as
the $\sigma$-function (up to a linear shift), equation (\ref{322}) coincides, up to the changing
the independent variable $s$ to the variable $t = \sqrt{s}$,  with  equation (C.29) of \cite{JM}
{\footnote{When comparing equation (\ref{322})  with equation (C.29) of \cite{JM}
one has to also take into account that, if we denote the $\sigma$-function
of \cite{JM}, for which equation (C.29) is written, as $\sigma_{JM}$ then the exact
relation with our $\sigma$ function is given by the equation,
$\sigma_{JM}(t) = 2\sigma(t^2).$}} .

\setcounter{equation}{0}
\section{Discrete $\sigma-$ form}

We may anticipate due to the recurrence relations and non-linear difference
equations, (2.16) and (2.17) that,
for a fixed $t,$ $H_n,$ and $H_{n\pm 1}$
would satisfy a discrete analog of (3.24). It turns out that in this instance $a_n$ has a simpler
expression in term of $H_{n\pm 1}.$ We shall obtain an expression for $b_n$ in terms of $H_{n}$ and
$H_{n\pm 1}.$ From
$$
H_n=-\sum_{j=0}^{n-1}a_j,
$$
we have
\bea
a_n=H_n-H_{n+1},
\eea
and
\bea
a_{n}+a_{n-1}=H_{n-1}-H_{n+1}=:\delta^2H_n.\nonumber
\eea
Multiply the above by $\bt_n$ we find,
\bea
\bt_n(a_n+a_{n-1})&=&\bt_n\dth\\
&=&ns-(2n+\al)b_n,
\eea
where the last equation follows from (2.12).
 Now (2.11) becomes
$$
\bt_n=n(n+\al)+b_n-H_n.
$$
Substituting the above into
(4.2) produces a linear equation in $b_n$ whose solution is
\bea
b_n=\frac{ns+\dth[H_n-n(n+\al)]}{2n+\al+\dth}.
\eea
Hence the auxiliary quantities $a_n$ and $b_n,$ and $\bt_n$ are now expressed
in terms of $H_n,$ and $H_{n\pm 1}.$ Substituting these into (2.13)
give rise to the discrete $\sigma-$form stated in the next Theorem.
\begin{theorem}
If
$$
H_n:=s\frac{d}{ds}\ln D_n,
$$
then
$$
\{[H_n-n(n+\al)]\dth+ns\}\{[H_n-n(n+\al)-s]\dth-(n+\al)s\}
$$
\begin{equation}\label{discrete}
=(2n+\al+\dth)\{ns+(2n+\al)[n(n+\al)-H_n]\}(H_n-H_{n+1})(H_{n-1}-H_n).
\end{equation}
\end{theorem}
We also have the discrete analog of (3.28) and (3.29),
\bea
\al_n&=&2n+1+\al+H_n-H_{n+1}\\
\bt_n&=&n(n+\al)+\frac{ns-n(n+\al)\dth-(2n+\al)H_n}{2n+\al+\dth}.
\eea
Now the obvious equalities, $(4.6)=(3.28)$ and $(4.7)=(3.29),$ imply two further differential-difference
equations which $H_n$ must satisfy,
\bea
H_n-H_{n+1}&=&\frac{2s(H_n')^2-sH_n'}{sH_n''+n-(2n+\al)H_n'}\\
\frac{ns-n(n+\al)\dth-(2n+\al)H_n}{2n+\al+\dth}&=&sH_n'-H_n.
\eea
\vskip .3in
\noindent
{\bf Remark IV.\:}  From the point of view of general Jimbo-Miwa
isomonodromy theory of Painlev\'e equations,
which we will be discussing in sections 5 and 6, equation (\ref{discrete})
should be related to  the B\"acklund - Schlesinger transformations of the
$\tau$ - function. However, we failed to identified equation (\ref{discrete})
with any of the difference equations for the $\tau$ - function discussed
in \cite{JM} and describing the possible Schlesinger transformations.
Since it is written for the logarithmic derivative of the $\tau$ - function,
and not for the $\tau$ - function itself as in \cite{JM}, equation (\ref{discrete})
might be in fact of a different nature than the ones considered in \cite{JM}.
We also want to mention that, equation (\ref{discrete}) is  an integrable discrete equation -
its Lax pair is formed by the first and the third equations of the triple (\ref{lax}) of section 5,
and the equation itself, as it has already been noticed, represents a
B\"acklund - Schlesinger transformation of the third Painlev\'e equation.
Therefore, we expect this equation to be equivalent to one of the known discrete Painlev\'e
equations, which we have not yet identify. Apparently this identification
is not quite  straightforward. One of the referees has suggested that (4.5) could be
a composition of the basic Schlesinger transformation $T_1$ and $T_2$
found in \cite{F+W2} (see Proposition (4.6) of \cite{F+W2}).

\setcounter{equation}{0}
\section{An alternative derivation of the Painlev\'e III equation.}

In this section we present an alternative derivation of the third
Painlev\'e equation (\ref{p3}) for the quantity $a_{n}(s).$ This
derivation is based on the Riemann-Hilbert point of view \cite{FIK4}, \cite{FIK5}  on
orthogonal polynomials and makes use of
the  general Jimbo-Miwa-Ueno theory of isomonodromy deformations.
This in turn allows us to place some of the key identities
of the preceding sections into a general framework of integrable systems.

The  Riemann-Hilbert problem for the orthogonal polynomials
at hand is the following

\begin{itemize}
\item $Y : \mathbb C \setminus \mathbb R_{+} \to
    \mathbb C^{2\times 2}$ is analytic.
\item $Y_+(x) = Y_-(x)
    \begin{pmatrix} 1 & x^{\alpha} e^{-x-s/x} \\ 0 & 1 \end{pmatrix}$
    for $x \in \mathbb R_{+} \setminus \{0\}$, with $\mathbb R_{+}$ oriented from left to right.
\item $Y(z) = (I + O(1/z)) \begin{pmatrix} z^n & 0 \\ 0 & z^{-n} \end{pmatrix}$
    as $z \to \infty$.
\item $Y(z) = O\left(1\right)$ as $z \to 0$.
\end{itemize}
Here $Y_{\pm}(z)$ denote the non-tangential limiting values of $Y(z)$
on $\mathbb R_{+}$ taken (in the usual point-wise sense) from the $\pm$ - side.
The Riemann-Hilbert problem has the unique solution expressed
in terms of the orthogonal polynomials $P_{n}(z)$,
\begin{equation} \label{RHsolution}
Y(z) = \begin{pmatrix}
    \ds P_{n}(z) &
    \ds \frac{1}{2\pi i}\,  \int_{\mathbb R_{+}} \frac{P_{n}(x) x^{\alpha} e^{-x -s/x}}{x-z} dx \\[10pt]
    -\frac{2\pi i}{h_{n-1}}\,P_{n-1}(z) &
    \ds -\frac{1}{h_{n-1}} \int_{\mathbb R_{+}} \frac{P_{n-1}(x) x^{\alpha} e^{-x-s/x}}{x-z} dx
  \end{pmatrix}.
\end{equation}
We also note that an immediate consequence of the uni-modularity of
the jump matrix of the above Riemann-Hilbert  problem
is the identity{\footnote{A more conventional derivation of this identity
is based on the use of the basic three-term recurrence  equation (\ref{rec})
(see e.g.\cite{CI})}},
\begin{equation}\label{det}
\det Y(z) \equiv 1.
\end{equation}

Equation (\ref{RHsolution}) implies, in particular, that the asymptotic behavior of the
function  $Y(z)$ at $z = \infty$ and $z=0$ can be specified as the following full
asymptotic series,
\begin{equation}\label{serinfty}
Y(z) \sim \left(I + \sum_{k=1}^{\infty}\frac{Y_{-k}}{z^{k}}\right)z^{n\sigma_{3}},
\quad z \to \infty,
\end{equation}
and
\begin{equation}\label{serzero}
Y(z) \sim Q\left(I + \sum_{k=1}^{\infty}Y_{k}z^{k}\right),
\quad z \to 0,
\end{equation}
where $\sigma_{3}$  denotes, as usual,  the third Pauli matrix
$$
\sigma_{3} =
 \begin{pmatrix}
    \ds 1 &0 \\[10pt]
    0&-1 \end{pmatrix}.
$$
Moreover, the matrix coefficients $Y_{\pm k}$ and $Q$ of these
series are the smooth functions
of $n$ and $s$ (and of course of $\alpha$), and  they all can be easily  expressed in terms
of the fundamental objects
associated with the orthogonal polynomials $P_{n}(z)$, i.e. in terms of the functions
$\mathsf{p}_k(n,s), h_{n}(s),$ and the negative
moments of the polynomials $P_{n}(z)$. Indeed, by a straightforward
calculation we have from (\ref{RHsolution}) the following expressions  for the coefficient
$Y_{-1}$ and for the matrix multiplier $Q \equiv Y(0)$.
\begin{equation}\label{Y-1}
Y_{-1} =
\begin{pmatrix}
    \mathsf{p}_1(n) &-\frac{h_{n}}{2\pi i}\\[10pt]
    -\frac{2\pi i}{h_{n-1}}&- \mathsf{p}_1(n)\end{pmatrix},
\end{equation}
\begin{equation}\label{Q}
Q =
\begin{pmatrix}
    1 &p_{n}\\[10pt]
    -q_{n}&1-p_{n}q_{n}\end{pmatrix}
    P^{\sigma_{3}}_{n}(0),
\end{equation}
where
\begin{equation}\label{pq}
p_{n} = \frac{P_{n}(0)}{2\pi i}\,  \int_{\mathbb R_{+}} \frac{P_{n}(x) x^{\alpha} e^{-x -s/x}}{x} dx,
\quad
q_{n} = \frac{2\pi i}{h_{n-1}} \frac{P_{n-1}(0)}{P_{n}(0)},
\end{equation}
and we have taken into account the determinant identity (\ref{det}).

We are now going to write down the triple of the differential and difference
equations for the function $Y_{n}(z,s)$ following the standard procedure
of the theory of integrable systems (see \cite{FT},  \cite{JMU}, \cite{JM}; see also \cite{FIK4},
\cite{ITW}, \cite{CI}, and \cite{FIKN}).

Put
\begin{equation}\label{Psidef}
\Psi(z) \equiv \Psi(z,n,s,\alpha)
:=Y(z)e^{-\frac{1}{2}\left(z + \frac{s}{z}\right)\sigma_{3}}z^{\frac{\alpha}{2}\sigma_{3}},
\end{equation}
where the branch of the function $z^{\frac{\alpha}{2}}$ we define by the
condition, $-\pi < \arg z < \pi $.
The Riemann-Hilbert relations in terms of the function $\Psi(z)$ reads
as follows,
\begin{itemize}
\item[{\bf{(i)}}] $\Psi : \mathbb C \setminus \mathbb R \to
    \mathbb C^{2\times 2}$ is analytic.
\item[{\bf{(ii)}}] $\Psi_+(x) = \Psi_-(x)
    \begin{pmatrix} 1 & 1 \\ 0 & 1 \end{pmatrix}$
    for $x \in \mathbb R$ and $ x > 0$.
\item[{\bf{(iii)}}] $\Psi_+(x) = \Psi_-(x)e^{\pi i\alpha \sigma_{3}}$
    for $x \in \mathbb R$ and $ x < 0$.
\item[{\bf{(iv)}}]$\Psi(z) \sim \left(I + \sum_{k=1}^{\infty}\frac{\Psi_{-k}}{z^{k}}\right)
z^{\left(n+\frac{\alpha}{2}\right)\sigma_{3}}e^{-\frac{z}{2}\sigma_{3}}$
as $z \to \infty$.
\item[{\bf{(v)}}]$\Psi(z) \sim Q\left(I + \sum_{k=1}^{\infty}\Psi_{k}z^{k}\right)
z^{\frac{\alpha}{2}\sigma_{3}}e^{-\frac{s}{2z}\sigma_{3}}$
as $z \to 0$,
\end{itemize}
where the real line  $\mathbb R$ is oriented as usual from left to right. The
coefficients $\Psi_{\pm k}$ of the asymptotic series  are easily evaluated
via combinations of  the coefficients $Y_{\pm k}$  and the coefficients
of the expansions of the exponential function $e^{-\frac{s}{2z}}$ and
$e^{-\frac{z}{2}}$ near $z =\infty$ and $z =0$, respectively. In particular, we have that,
\begin{equation}\label{Psi-1}
\Psi_{-1} =
\begin{pmatrix}
    \mathsf{p}_1(n) -\frac{s}{2} &-\frac{h_{n}}{2\pi i}\\[10pt]
    -\frac{2\pi i}{h_{n-1}}&- \mathsf{p}_1(n) + \frac{s}{2}\end{pmatrix}.
\end{equation}
The important feature of the $\Psi$-RH problem is that the jump matrices
of the jump relations {\bf(ii)} and {\bf(iii)}
do not depend on $z$,  $s$, and $n$. Therefore, by standard
arguments based on the Liouville theorem (cf. \cite{FT}, \cite{FIK4}), we conclude
that the logarithmic derivatives,
\begin{equation}\label{AVU}
A(z) := \frac{\partial \Psi(z)}{\partial z}\Psi^{-1}(z), \quad
B(z) := \frac{\partial \Psi(z)}{\partial s}\Psi^{-1}(z), \quad\mbox{and}\quad
U(z) := \Psi(z,n+1)\Psi^{-1}(z,n),
\end{equation}
are rational functions of $z$. Using the asymptotic expansions {\bf(iv)} and {\bf(v)},
we can evaluate the respective principal parts at the poles
at the points $z=0$ and $z=\infty$ and arrive, taking into account
(\ref{Psi-1}) and (\ref{Q}), at the following explicit
formulae for the function $A(z)$, $B(z)$, and $U(z)$.
\begin{equation}\label{A}
A(z) = -\frac{1}{2}\sigma_{3} + \frac{A_{1}}{z} +  \frac{A_{2}}{z^2},
\end{equation}
\begin{equation}\label{V}
B(z) = -\frac{A_{2}}{sz},
\end{equation}
and
\begin{equation}\label{U}
U(z) = z\begin{pmatrix}1&0\\[10pt]0&0\end{pmatrix} + U_{0},
\end{equation}
where{\footnote{Notation $[M_1, M_2]$ means the usual commutator of two
matrices, $[M_1, M_2] = M_{1}M_{2} - M_{2}M_{1}$.}}
$$
A_{1} = \frac{1}{2}[\sigma_{3}, \Psi_{-1}] + \left(n+\frac{\alpha}{2}\right)\sigma_{3}
$$
\begin{equation}\label{A1}
=
\begin{pmatrix}
  n+ \frac{\alpha}{2} &-\frac{h_{n}}{2\pi i}\\[10pt]
    \frac{2\pi i}{h_{n-1}}& - n - \frac{\alpha}{2}\end{pmatrix} ,
\end{equation}
\vskip .2in

\begin{equation}\label{A2}
A_{2} = \frac{s}{2}Q\sigma_{3}Q^{-1}
=\frac{s}{2}\begin{pmatrix}
  1 -2p_{n}q_{n} &-2p_{n}\\[10pt]
    2q_{n}\left(p_{n}q_{n} - 1\right)&2p_{n}q_{n}-1 \end{pmatrix} ,
\end{equation}
and
$$
U_{0} = \Psi_{-1}(n+1)\begin{pmatrix}1&0\\[10pt]0&0\end{pmatrix}
-\begin{pmatrix}1&0\\[10pt]0&0\end{pmatrix}\Psi_{-1}(n)
$$
\begin{equation}\label{U0}
=
\begin{pmatrix}
 \mathsf{p}_1(n+1)-\mathsf{p}_1(n)  &\frac{h_{n}}{2\pi i} \\[10pt]
  -\frac{2\pi i}{h_{n}}   &0 \end{pmatrix}
  \equiv
\begin{pmatrix}
 -\alpha_{n}  &\frac{h_{n}}{2\pi i} \\[10pt]
 - \frac{2\pi i}{h_{n}}   &0 \end{pmatrix}
\end{equation}

According to the standard methodology (cf. \cite{JMU}, \cite{JM}, \cite{FIK4}),
relations (\ref{AVU}) should be now re-interprited as a
system of linear differential-difference equations,
\begin{equation}\label{lax}
 \left\{ \begin{array}{ll}
    \frac{\partial \Psi(z)}{\partial z} = A(z)\Psi(z) \\\\
    \frac{\partial \Psi(z)}{\partial s} = B(z)\Psi(z) \\\\
    \Psi(z,n+1) = U(z)\Psi(z,n),
\end{array}\right.
\end{equation}
which we call the Lax triple.

The compatibility conditions of this system, i.e. the equations,
\begin{equation}\label{comp1}
\frac{\partial A(z)}{\partial s} - \frac{\partial B(z)}{\partial z}
=[B(z),A(z)] \quad \Bigl(\Psi_{zs} = \Psi_{sz}\Bigr)
\end{equation}
\begin{equation}\label{comp2}
\frac{\partial U(z)}{\partial s}
= B(z,n+1)U(z) - U(z)B(z,n), \quad \Bigl(
(\Psi(z,n+1))_{s} = \Psi_{s}(z,n+1)\Bigr)
\end{equation}
and
\begin{equation}\label{comp3}
\frac{\partial U(z)}{\partial z}
= A(z,n+1)U(z) - U(z)A(z,n), \quad  \Bigl(
(\Psi(z,n+1))_{z} = \Psi_{z}(z,n+1)\Bigr)
\end{equation}
yield the  Painlev\'e type  (equ. (\ref{comp1})),
the Toda type  (equ. (\ref{comp2})) and the discrete
Painlev\'e or Freud type (equ. (\ref{comp3})) equations,
respectively,
for a proper combination of functions $h_{n}$, $p_{n}$ and
$q_{n}$. Moreover, using the Jimbo-Miwa list of the Lax pairs
for  Painlev\'e equations \cite{JM} and noticing that the ``master equation'',
i.e. the first equation of
system (\ref{lax}) is a $2\times 2$ system with two
irregular singular points of Poincare rank 1, one concludes
that the relevant Painlev\'e equation is, in fact, the third Painlev\'e
equation. In order to make a precise statement, i.e. to point out the
exact combination of  functions $h_{n}$, $p_{n}$ and
$q_{n}$ which makes the solution of  Painlev\'e III equation,
we  only need to perform
a simple scaling transformation of system (\ref{lax})
which would bring it to the normal form   of  \cite{JM}. To this end,
we introduce the new independent variables,
\begin{equation}\label{tlambda}
\lambda:=s^{-1/2}z, \quad \mbox{and}\quad t := \sqrt{s},
\end{equation}
and pass from the function $\Psi(z,s)$ to the function $\Phi(\lambda,t)$
defined by the equation,
\begin{equation}\label{defPhi}
\Phi(\lambda, t):= t^{-\left(n+\frac{\alpha}{2}\right)\sigma_{3}}\Psi(t\lambda,t^2).
\end{equation}
We notice that in  terms of the function $\Phi(\lambda, t)$ the asymptotic
relations {\bf {(iv)}} and {\bf {(v)}} transform into the relations,
\begin{equation}\label{Phiinfty}
\Phi(\lambda) \sim \left(I + \sum_{k=1}^{\infty}\frac{\Phi_{-k}}{\lambda^{k}}\right)
\lambda^{\left(n+\frac{\alpha}{2}\right)\sigma_{3}}e^{-\frac{t\lambda}{2}\sigma_{3}},
\quad \lambda \to \infty,
\end{equation}
and
\begin{equation}\label{Phizero}
\Phi(\lambda) \sim R\left(I + \sum_{k=1}^{\infty}\Phi_{k}\lambda^{k}\right)
\lambda^{\frac{\alpha}{2}\sigma_{3}}e^{-\frac{t}{2\lambda}\sigma_{3}},
\quad \lambda \to 0,
\end{equation}
with the new coefficients connected  to the old ones by the equations,
$$
R = t^{-\left(n+\frac{\alpha}{2}\right)\sigma_{3}}Qt^{\frac{\alpha}{2}\sigma_{3}},
\quad \Phi_{-k} = t^{-k}t^{-\left(n+\frac{\alpha}{2}\right)\sigma_{3}}\Psi_{-k}
t^{\left(n+\frac{\alpha}{2}\right)\sigma_{3}},
\quad \Phi_{k} = t^{k}t^{-\frac{\alpha}{2}\sigma_{3}}\Psi_{k}
t^{\frac{\alpha}{2}\sigma_{3}}.
$$
Simultaneously, the first two equations of
system (\ref{lax}) transform into the {\it Jimbo-Miwa Lax pair}
for the third Painlev\'e equation,
\begin{equation}\label{laxJM}
 \left\{ \begin{array}{ll}
    \frac{\partial \Phi(\lambda)}{\partial \lambda} =
    \left(-\frac{t}{2}\sigma_{3} +\frac{A_{-1}}{\lambda}
    + \frac{A_{-2}}{\lambda^2}\right)\Phi(\lambda)
    \equiv A_{JM}(\lambda)\Phi(\lambda) \\\\
  \frac{\partial \Phi(\lambda)}{\partial t} =
    \left(-\frac{\lambda}{2}\sigma_{3} +B_{0}
    + \frac{B_{-1}}{\lambda}\right)\Phi(\lambda)
    \equiv B_{JM}(\lambda)\Phi(\lambda).
\end{array} \right.
\end{equation}
Here, the matrix coefficients $A_{-1}$, $A_{-2}$, $B_{0}$, and $B_{-1}$
are given by the equations,
\begin{equation}\label{A-1}
A_{-1} = t^{-\left(n+\frac{\alpha}{2}\right)\sigma_{3}}
A_{1}t^{\left(n+\frac{\alpha}{2}\right)\sigma_{3}}
=
\begin{pmatrix}
  -\frac{\theta_{\infty}}{2} &u\\[10pt]
    v&\frac{\theta_{\infty}}{2} \end{pmatrix},
\end{equation}
\vskip .2in
\begin{equation}\label{A-2}
A_{-2} =\frac{1}{t} t^{-\left(n+\frac{\alpha}{2}\right)\sigma_{3}}
A_{2}t^{\left(n+\frac{\alpha}{2}\right)\sigma_{3}}
=
\begin{pmatrix}
  \zeta + \frac{t}{2} &-w\zeta\\[10pt]
   \frac{\zeta + t}{w}&- \zeta - \frac{t}{2} \end{pmatrix},
\end{equation}
\vskip .2in
\begin{equation}\label{B0}
B_{0} = \frac{1}{t}A_{-1} - \frac{1}{t} \left(n+\frac{\alpha}{2}\right)\sigma_{3}
= \frac{1}{t}\begin{pmatrix}
  0&u\\[10pt]
    v&0\end{pmatrix},
\end{equation}
\vskip .2in
\begin{equation}\label{B-1}
B_{-1} =-\frac{1}{t}A_{-2}
=
-\frac{1}{t}\begin{pmatrix}
  \zeta + \frac{t}{2} &-w\zeta\\[10pt]
   \frac{\zeta + t}{w}&- \zeta - \frac{t}{2} \end{pmatrix},
\end{equation}
where
\begin{equation}\label{thetainfty}
\theta_{\infty}= -\alpha -2n,
\end{equation}
and the new scalar functional parameters $u$, $v$, $\zeta$, and $w$ are defined in terms
of the original functions $h_{n}$, $h_{n-1}$, $p_{n}$, and $q_{n}$ via the formulae,
\begin{equation}\label{uv}
u=-\frac{h_{n}}{2\pi i}t^{-2n-\alpha}, \quad
v=\frac{2\pi i}{h_{n-1}}t^{2n+\alpha},
\end{equation}
\begin{equation}\label{zetaw}
\zeta = -tp_{n}q_{n},\quad w =-\frac{1}{q_{n}}t^{-2n-\alpha}.
\end{equation}

With the $u-w$ notations, the Lax pair (\ref{laxJM}) matches, up to
the replacement $t \to -t$ and the use of the letter $z$ instead of the letter $\zeta$,
the Lax pair presented on page 439 of \cite{JM}, and hence we can use the
general results of  \cite{JM}.

\begin{theorem}[\cite{JM}] \label{th4}
Consider the over-determined linear system (\ref{laxJM})
with the matrix coefficients defined by the equations
\begin{equation}\label{A-12}
A_{-1} =
\begin{pmatrix}
  -\frac{\theta_{\infty}}{2} &u\\[10pt]
    v&\frac{\theta_{\infty}}{2} \end{pmatrix},\quad
A_{-2} =
\begin{pmatrix}
  \zeta + \frac{t}{2} &-w\zeta\\[10pt]
   \frac{\zeta + t}{w}&- \zeta - \frac{t}{2} \end{pmatrix},
\end{equation}
\vskip .2in
\begin{equation}\label{B0-1}
B_{0} = \frac{1}{t}\begin{pmatrix}
  0&u\\[10pt]
    v&0\end{pmatrix}, \quad
B_{-1} =
-\frac{1}{t}\begin{pmatrix}
  \zeta + \frac{t}{2} &-w\zeta\\[10pt]
   \frac{\zeta + t}{w}&- \zeta - \frac{t}{2} \end{pmatrix} ,
\end{equation}
i.e., by the right hand sides of the last equalities in the formulae
(\ref{A-1}) - (\ref{B-1}) not necessarily    assuming any specific,  ``orthogonal polynomial''  choice
of the parameters $\theta_{\infty}$,
$u, v, \zeta$ and $w$. Then the following is true.
\begin{enumerate}
\item The compatibility condition of system
(\ref{laxJM}), i.e. the matrix equation,
\begin{equation}\label{lax2}
\frac{\partial A_{JM}}{\partial t} - \frac{\partial B_{JM}}{\partial \lambda} =
[B_{JM},A_{JM}],
\end{equation}
is equivalent to the following set of scalar equations,
\begin{equation}\label{lax3}
t\frac{du}{dt} = \theta_{\infty}u +2t\zeta w,
\end{equation}
\begin{equation}\label{lax4}
t\frac{dv}{dt} = -\theta_{\infty}v +\frac{2t}{w}(\zeta + t),
\end{equation}
\begin{equation}\label{lax5}
t\frac{d\zeta}{dt} = 2\zeta wv + \zeta + \frac{2u(\zeta + t)}{w},
\end{equation}
\begin{equation}\label{lax6}
t\frac{d\ln w}{dt} = \frac{2u}{w} -2wv - \theta_{\infty},
\end{equation}
with the quantity,
\begin{equation}\label{theta0def}
\theta_{0}: = -\frac{\theta_{\infty}}{t}(2\zeta +t) + \frac{2u(\zeta + t)}{tw}
-\frac{2\zeta}{t}wv,
\end{equation}
being  the first integral of system (\ref{lax3}) - (\ref{lax6}).
\item The quantity  $\theta_{0}$
is very similar to the parameter $\theta_{\infty}$. Indeed, they describe
the {\it formal monodromy} at the relevant irregular points;
$\theta_{\infty}$ -- at $\lambda = \infty$ and $\theta_{0}$ -- at $\lambda = 0.$
This means that they appear as the branching exponents in
the following formal matrix solutions{\footnote{ All the coefficients
$\Phi^{(\infty)}_{k}$ of series (\ref{PhiinftyJM}) are uniquely defined as
rational functions of $u, v, \zeta, w,$ and $t$ via simple recurrence relations \cite{JMU}.
The matrix factor $R^{(0)}$ is defined, up to a right diagonal multiplier,
$R^{(0)} \to R^{(0)}\Lambda$,
by the equations, $\frac{t}{2}R^{(0)}\sigma_{3}\left(R^{(0)}\right)^{-1} =
A_{-2},\quad \det R^{(0)} =1$. The coefficients $\Phi^{(0)}_{k}$ of series (\ref{PhizeroJM})
are defined  as rational functions of $u, v, \zeta, w,$ and $t$ up to the conjugation
$\Phi^{(0)}_{k} \to {\Lambda}^{-1} \Phi^{(0)}_{k} \Lambda$. We will
discuss this issue in more detail in the next section.}}    of the  Lax pair (\ref{laxJM})
at the points $\lambda =\infty$ and  $\lambda = 0$.
\begin{equation}\label{PhiinftyJM}
\Phi^{(\infty)}_{\mbox{formal}}(\lambda) =
\left(I + \sum_{k=1}^{\infty}\frac{\Phi^{(\infty)}_{k}}{\lambda^{k}}\right)
\lambda^{-\frac{\theta_{\infty}}{2}\sigma_{3}}e^{-\frac{t\lambda}{2}\sigma_{3}},
\end{equation}
and
\begin{equation}\label{PhizeroJM}
\Phi^{(0)}_{\mbox{formal}}(\lambda) =
R^{(0)}\left(I + \sum_{k=1}^{\infty}\Phi^{(0)}_{k}\lambda^{k}\right)
\lambda^{\frac{\theta_{0}}{2}\sigma_{3}}e^{-\frac{t}{2\lambda}\sigma_{3}}.
\end{equation}
\item The function
 \begin{equation}\label{p3JM}
y := -\frac{u}{\zeta w},
\end{equation}
satisfies   the third Painlev\'e equation,
\begin{equation}\label{p3JM}
\frac{d^2y}{dt^2} = \frac{1}{y}\left(\frac{dy}{dt}\right)^2
- \frac{1}{t}\frac{dy}{dt} - \frac{1}{t}\Bigl(4\theta_{0}y^2 + 4(1-\theta_{\infty})\Bigr)
+4y^3 -\frac{4}{y}.
\end{equation}
In the standard  notations \cite{ince}, this is the Painlev\'e III equation,
$P_{III}(-4\theta_{0}, -4(1-\theta_{\infty}), 4, - 4)$.
\end{enumerate}
\end{theorem}

Comparing  the general formal expansions  (\ref{PhiinftyJM}), (\ref{PhizeroJM})
with the asymptotic series (\ref{Phiinfty}), (\ref{Phizero}), we conclude that
in our case the parameters $\theta_{\infty}$ and $\theta_{0}$ assume the
values (see also (\ref{thetainfty})),
\begin{equation}\label{thetainftyzero}
\theta_{\infty}= -\alpha -2n,\quad \mbox{and}\quad \theta_{0} = \alpha.
\end{equation}
Simultaneously, from (\ref{uv})-(\ref{zetaw}) we see that in our case,
\begin{equation}\label{your}
y = \frac{h_{n}}{2\pi itp_{n}}.
\end{equation}
In other words, relation (\ref{your}) provides us with the combination
of the functions $h_{n}, p_{n}, q_{n}$ which satisfies the third Painlev\'e
equations and which we have been looking for, while relations
(\ref{thetainftyzero}) specifies  the parameters of the Painlev\'e
equation. We conclude then that the function $y$ defined in (\ref{your})
satisfies the following Painlev\'e III equation,
\begin{equation}\label{p3y}
\frac{d^2y}{dt^2} = \frac{1}{y}\left(\frac{dy}{dt}\right)^2
- \frac{1}{t}\frac{dy}{dt} - \frac{1}{t}\Bigl(4\alpha y^2 + 4( 2n +1 + \alpha)\Bigr)
+4y^3 -\frac{4}{y}.
\end{equation}
 What is still left for us is to establish the connection
between the function $y(t)$ defined in (\ref{your}) and the function
$a_{n}(s)$ defined in (\ref{andef}). This is easy; indeed, we have that
$$
 \int_{0}^{\infty}\frac{P_n^2}{x}wdx =
 \int_{0}^{\infty}P(x)\left(x^{n-1} +   \mathsf{p}_1(n)x^{n-2}
 +...+\frac{P_{n}(0)}{x}\right)wdx
 = P_{n}(0)\int_{0}^{\infty}\frac{P_n}{x}wdx,
 $$
 and hence (see (\ref{pq})),
\begin{equation}\label{pnan}
p_{n} = \frac{h_{n}}{2\pi i}\frac{a_{n}}{s},
\end{equation}
which in turn implies that
\begin{equation}\label{yan}
y = \frac{t}{a_{n}} \equiv \frac{X_{n}}{t}.
\end{equation}
It is an elementary exercise to check that the
substitution,
$$
y =  \frac{t}{a_{n}}, \quad t = \sqrt{s},
$$
transforms equation (\ref{p3y}) into equation (\ref{p3}) for the quantity $a_{n}(s)$.

We conclude this section by revealing  the connections
of some of the key identities established in the previous sections with
the general constructions of the isomonodromy theory
of Painlev\'e equations discussed above.

We first notice that, similar to the derivation of equation (\ref{pnan}), the following
\begin{equation}\label{qnan}
q_{n} = \frac{2\pi i}{h_{n}}\frac{b_{n}}{a_{n}}.
\end{equation}
Therefore all the matrix coefficients of the Lax  triple  (\ref{lax})
can be expressed in terms of the functions $\alpha_n, h_{n}, a_{n},$
and $b_n$. Actually, we need only to re-write the coefficient
$A_2$ from (\ref{A2}),
\begin{equation}\label{A2new}
A_{2}
=\frac{s}{2}\begin{pmatrix}
  1 -\frac{2b_{n}}{s} &-\frac{h_{n}}{\pi i}\frac{a_{n}}{s}\\[10pt]
    \frac{4\pi i}{h_{n}}\frac{b_{n}}{a_{n}}\left(\frac{b_{n}}{s} - 1\right)&
   \frac{2b_{n}}{s} -1 \end{pmatrix}.
\end{equation}
Substituting this representation together with the similar formulae
for $A_1$  and $U_{0}$ (see (\ref{A1}) and (\ref{U0})) into the
compatibility  equation (\ref{comp3}) yields the set of the scalar
difference equations (\ref{29}) - (\ref{213}) of section 2. Similar operation with
the compatibility equation (\ref{comp2}) results in the
differential-difference equations  (\ref{31}) - (\ref{36}) of section 3.
The scalar form of compatibility equation  (\ref{comp1}) we have already discussed
in detail in this section. As we have seen, this equation is equivalent
to the set of scalar equations (\ref{lax3}) - (\ref{lax6}) which,
in terms of the functions $h_{n}(s),h_{n-1}(s), a_{n}(s), b_{n}(s)$, transforms
to the system,
\begin{equation}\label{lax3h}
s\frac{dh_{n}}{ds} = -h_{n}a_{n},
\end{equation}
\begin{equation}\label{lax4h}
s\frac{dh_{n-1}}{ds} = \frac{h_{n-1}^2}{h_{n}}\frac{b_{n}}{a_{n}}(s-b_{n}),
\end{equation}
\begin{equation}\label{lax5b}
s\frac{db_{n}}{ds} = b_{n} -\frac{b_{n}}{a_{n}}(s-b_{n}) - \frac{h_{n}}{h_{n-1}}a_{n}
\end{equation}
\begin{equation}\label{lax6a}
s\frac{da_{n}}{ds} = 2b_{n} + (2n+1 +\alpha + a_{n})a_{n} -s,
\end{equation}
which is equivalent to the system of equations derived in Lemma 5 of section 3.

We want  to highlight the theoretical meaning of the important
identity (\ref{214}). It is, in fact, the formal monodromy identity (\ref{theta0def})
written in terms of the functions $h_{n}(s),h_{n-1}(s), a_{n}(s), b_{n}(s)$.

It is also worth mentioning, that the ladder equations (\ref{21}) and (\ref{22})
are just the first column of the  master, $z$ - equation of the Lax triple (\ref{lax}).
The  first column of the second equation of the triple (\ref{lax})  yields
the ladder operators in $s$,
\begin{equation}\label{ladders1}
\Bigl(zs\frac{d}{ds} - b_{n}\Bigr)P_{n}(z) = -\beta_{n}a_{n}P_{n-1}(z),
\end{equation}
\begin{equation}\label{ladders2}
\Bigl(zs\frac{d}{ds} - b_{n-1} -\alpha_{n-1}a_{n-1}
+za_{n-1}\Bigr)P_{n-1}(z) = a_{n-1}P_{n}(z).
\end{equation}
When deriving these equations we made use of the relations (\ref{210}) and (\ref{213}).
The equations (\ref{ladders1}) and (\ref{ladders2}) of course can be derived using
the orthogonality relations, bypassing the isomonodromy theory.

Finally, the third equation in (\ref{lax}) is equivalent to the basic recurrence
relations (\ref{rec}) for the polynomials  $P_{n}$. It should be pointed out, that for the general
solutions of the Painlev\'e III equations, only the  first two equations, which are
equivalent under scaling to the Lax pair (\ref{laxJM}), take place. From the point of view of the general
isomonodromy theory of Painlev\'e equations, the  third equation of
the triple (\ref{lax}) describes the B\"acklund-Schlesinger transformation
of the third Painlev\'e equation (see \cite{JMU}; see also Chapter 6  of \cite{FIKN}).

\setcounter{equation}{0}
\section{The Hankel determinant as the isomonodromy $\tau$ - function.}

Let us remind the reader the Jimbo-Miwa definition of the $\tau$-function corresponding
to the third Painlev\'e equation (\ref{p3JM}).

Consider the formal series (\ref{PhiinftyJM}) and (\ref{PhizeroJM}) of theorem \ref{th4}.
In the footnote to this theorem we have
already mentioned that all the coefficients of these series can be evaluated
as rational functions of the parameters $u, v, \zeta, w,$ and $t$. For example,
the first  coefficients in each of he series  are given by the relations (see p. 440 of \cite{JM}),
\begin{equation}\label{Phi1infty}
\Phi^{(\infty)}_{1}
=\frac{1}{t}\begin{pmatrix}
  uv -t\zeta -\frac{t^2}{2} &u\\[10pt]
    -v&
    -uv +t\zeta +\frac{t^2}{2} \end{pmatrix} ,
\end{equation}
and
\begin{equation}\label{Phi1zero}
\Phi^{(0)}_{1}
=\frac{1}{t}\begin{pmatrix}
  \tilde{u}\tilde{v} -t\zeta -\frac{t^2}{2} &-\tilde{u}\\[10pt]
    \tilde{v}&
    -\tilde{u}\tilde{v} +t\zeta +\frac{t^2}{2} \end{pmatrix} ,
\end{equation}
where the parameters $\tilde{u}$ and $\tilde{v}$ are related to
the basic parameters $u$ and $v$ through the equation,
\begin{equation}\label{uvtilde}
\left(R^{(0)}\right)^{-1}
\begin{pmatrix}
 - \frac{\theta_{\infty}}{2} &u\\[10pt]
 v&\frac{\theta_{\infty}}{2} \end{pmatrix}R^{(0)}
 = \begin{pmatrix}
  \frac{\theta_{0}}{2} &\tilde{u}\\[10pt]
 \tilde{v}&-\frac{\theta_{0}}{2} \end{pmatrix},
 \end{equation}
which, in particular, means the identity,
\begin{equation}\label{uvtilde2}
\tilde{u}\tilde{v} = uv +\frac{\theta^{2}_{\infty} - \theta^{2}_{0}}{4}.
\end{equation}
Denote $\hat{Y}_{\infty}(\lambda)$ and $\hat{Y}_{0}(\lambda)$ the
series in the brackets of formulae  (\ref{PhiinftyJM}) and (\ref{PhizeroJM}),
i.e.
\begin{equation}\label{hatinfty}
\hat{Y}_{\infty}(\lambda) = \left(I + \sum_{k=1}^{\infty}\frac{\Phi^{(\infty)}_{k}}{\lambda^{k}}\right),
\end{equation}
and
\begin{equation}\label{hatzero}
\hat{Y}_{0}(\lambda) = \left(I + \sum_{k=1}^{\infty}\Phi^{(0)}_{k}\lambda^{k}\right).
\end{equation}
The Jimbo-Miwa-Ueno isomonodromy $\tau$-function \cite{JMU}, in the case of the
Lax pair (\ref{laxJM}) is defined by the formula,
\begin{equation}\label{jmutau}
d\ln \tau = -\mbox{Trace}\,\mbox{Res}_{\lambda=0}
\hat{Y}^{-1}_{0}(\lambda)\frac{\partial \hat{Y}_{0}}{\partial \lambda}(\lambda)
dT_{0}(\lambda)
-\mbox{Trace}\,\mbox{Res}_{\lambda=\infty}
\hat{Y}^{-1}_{\infty}(\lambda)\frac{\partial \hat{Y}_{\infty}}{\partial \lambda}(\lambda)
dT_{\infty}(\lambda),
\end{equation}
where
$$
dT_{0}(\lambda) = -\frac{1}{2\lambda}\sigma_{3}dt, \quad
\mbox{and}\quad dT_{\infty}(\lambda) = -\frac{\lambda}{2}\sigma_{3}dt.
$$
Substituting   (\ref{hatinfty}) and (\ref{hatzero}) into (\ref{jmutau}), we arrive
at the equation,
$$
d\ln \tau = \frac{1}{2}\mbox{Trace}\left(\Phi^{(0)}_{1}\sigma_{3}
+ \Phi^{(\infty)}_{1}\sigma_{3}\right)dt,
$$
which, taking into account (\ref{Phi1infty}), (\ref{Phi1zero}), and (\ref{uvtilde2}),
implies that,
\begin{equation}\label{tautau}
d\ln \tau(t) = H_{III}(u(t), v(t), \zeta(t);t)dt,
\end{equation}
where
\begin{equation}\label{HIII0}
H_{III}(u, v, \zeta; t)
= \frac{1}{t} \left(2uv  - 2t\zeta - t^2
- \frac{\theta^2_{0} - \theta^2_{\infty}}{4}\right).
\end{equation}

Let us now turn to the Hankel determnant $D_{n}(s)$ and
consider the quantity
$$
H_n:=\sps\ln D_n,
$$
which has played a central role in sections 3 and 4.
From (\ref{38}) and (\ref{215}) we have that
\begin{equation}\label{Hnour}
H_n = n(n+\alpha) + b_{n} - \frac{ns-(2n+\alpha)b_{n}}{a_{n}}
+\frac{b^2_{n} -sb_{n}}{a^2_{n}}.
\end{equation}
At the same time, in our case (see (\ref{thetainftyzero})),
$$
\theta_{0} = \alpha, \quad \theta_{\infty} = -\alpha  - 2n,
\quad uv = -\beta_{n}, \quad\mbox{and}\quad \zeta =  -\frac{b_{n}}{t},
$$
where for the last two formulae we have used (\ref{uv}), (\ref{zetaw}),
(\ref{pnan}) and (\ref{qnan}). Substituting these formulae into
(\ref{HIII0}) we obtain that for our special solution of
system (\ref{lax3}) - (\ref{lax6}) the function $H_{III}$ assumes the form,
\begin{equation}\label{HIII1}
H_{III} =  \frac{1}{t} \Bigl(-2\beta_{n} + 2b_{n}
-s +n(n+\alpha) \Bigr).
\end{equation}
Recalling now identity (\ref{214}) (which we remind the reader is the formal
monodromy identity (\ref{theta0def}) in disguise) we arrive,
after some simple algebra, at  the relation,
\begin{equation}\label{tau3}
H_{III} = \frac{2}{t} H_{n} - t - \frac{n(n+\alpha)}{t}.
\end{equation}
This equation, in its turn, implies that
$$
s\frac{d}{ds}\ln \tau = \frac{t}{2}\frac{d}{dt}\ln \tau
= \frac{t}{2}H_{III} = H_{n} - \frac{s}{2} - \frac{n(n+\alpha)}{2}
$$
\begin{equation}\label{tau4}
=s\frac{d}{ds}\ln D_{n} - \frac{s}{2} - \frac{n(n+\alpha)}{2},
\end{equation}
and hence we obtain the following relation between the
Hankel determinant $D_{n}$ and the $\tau$-function of the
third Painlev\'e equation (compare with  \cite{ITW} where
similar formula is derived for a class of Toeplitz determinants),
\begin{equation}\label{tau5}
D_{n}(s) = \mbox{const}\,\tau(s) e^{\frac{s}{2}}s^{\frac{n(n+\alpha)}{2}}.
\end{equation}

\noindent
{\bf Remark V.\:} The function $H_{III}(u,v,\zeta; t)$ actually depends on $\zeta$
and the product $uv$. The latter, can be expressed, with the help of the
formal monodromy relation   (\ref{theta0def}) in terms of $\zeta$, $t$, and $y$.
Indeed, we have that
$$
2uv = (\theta_{0} + \theta_{\infty})ty +2\zeta y \theta_{\infty} + 2y^2\zeta(\zeta + t),
$$
and expression (\ref{HIII0}) can be transformed to the following equation
defined $H_{III}$ as a function n $y$, $\zeta$ and $t$,
\begin{equation}\label{HIIIfinal}
H_{III}(y,\zeta; t) = \frac{1}{t}\left(2y^2\zeta^2 + 2(ty^2 +\theta_{\infty}y
-t)\zeta + (\theta_{0} + \theta_{\infty})ty -t^2 - \frac{\theta^2_{0} - \theta^2_{\infty}}{4}\right).
\end{equation}
This is the canonical (see again p. 440 of \cite{JM}) representation
of the logarithmic derivative of the $\tau$ - function for  Painlev\'e III
equation (\ref{p3JM}). The remarkable fact of the general Jimbo-Miwa-Ueno
theory is that the function $H_{III}(y,\zeta; t)$ is the Hamiltonian of the third Painlev\'e
equation.
\vskip .2in
\noindent
{\bf Remark VI.\:} In sections 2 and 3, an important role has been played by
equation (\ref{215}) which transforms a {\it non-local} object - the sum of
$a_{j}$ from $j =0$ to $j = n-1$, to a {\it local} expression, which involves only $a_{n}$ and
$b_{n}$. We can see now an intrinsic reason for that. Indeed, on one hand,
the sum mentioned is, by its very nature, the logarithmic derivative of the Hankel determinant.
On the other hand, the latter is a $\tau$ - function and hence its logarithmic derivative
{\it must} admit a local representation in view of the general formula
(\ref{jmutau}).
\vskip .2in
\noindent
{\bf Remark VII.\:} The isomondromy context for orthogonal polynomials
and, in particular, the interpretation of the Hankel determinants as  isomonodromy
$\tau$ - functions  have been well understood for some time, since the early nineties
works \cite {FIK4} and  \cite{magnus2}. For the most general semiclassical weights
this fact was established in the recent paper \cite{bertola} (see also \cite{beh})

\section*{Acknowledgements}

Yang Chen is supported  by EPSRC grant \# R27027.
Alexander Its was supported by EPSRC grants \# R27027
and \# EP/F014198/1
and by NSF grant \#DMS-0701768.

\end{document}